\documentclass[12pt]{article}

\usepackage{times}
\usepackage{bm}
\usepackage{graphicx}
\usepackage{amsfonts}
\usepackage{amsmath}
\usepackage{amssymb}
\usepackage{bm}
\usepackage{graphicx}

\usepackage[round]{natbib}
\usepackage{setspace,cite}
\begin{document}

\title{Effect of sampling on the estimation of drift parameter of continuous time AR(1) processes}
\author{Radhendushka Srivastava \footnote{Radhendushka Srivastava was a postdoctorate researcher supported by NSF-DMS 0808864, NSF-EAGER 1249316, a gift from Microsoft, a gift from Google, and Ping Li's salary recovery account.} \and Ping Li \\
{\centerline{\small Department of Statistical Science, Cornell University}} }
\maketitle

\begin{abstract}
We study the effect of stochastic sampling on the estimation of the drift parameter of continuous time AR(1) process. A natural distribution free moment estimator is considered for the drift based on stochastically observed time points. The effect of the constraint of the minimum separation between successive samples on the estimation of the drift is studied.
\end{abstract}

\section{Introduction}

Sampling is an integral component of the inference for continuous time processes. The stochastic diffusion equations are the important continuous time stochastic models. These diffusion equations are often used to model in economic studies \citep{SM,DG,Fan,TC}, in communication theory \citep{CV,STP,MT,HMR}, in sea surface temperature data \citep{TAA}, and several other discipline of science and engineering. The continuous time AR(1) process is the first order stochastic diffusion equation (also known as Ornstein-Uhlenbeck (OU) process). This is the simplest stochastic model used in sensor network \citep{HMR}, in sea surface temperature data \citep{TAA}, in financial time series \citep{TC}, in physics \citep{Chandrasekhar}, in meteorology \citep{Gringorten}, in population growth model \citep{Tuckwell} and in neurophysiological study \citep{ Stein,Linetsky} etc.

The continuous time AR(1) process is defined as the stationary solution of the stochastic differential equation
\begin{equation}\label{AReqn}
\frac{d}{dt}X_t=-\alpha X_t +\sigma W_t,
\end{equation}
where $W_t$ is the stationary white noise, $\alpha>0$ and $\sigma>0$, and the derivative operator is interpreted in weak sense. The process $X$ is completely specified by the drift parameter $\alpha$ and innovation noise $\sigma^2$.

The common sampling strategy of the continuous time process is the uniform and the stochastic point processes. When the OU process is observed at uniformly spaced intervals (say $\delta$) time point, then the sampled process constitutes discrete time autoregressive process
\begin{equation}\label{AReqn1}
X_{t\delta}=\phi X_{(t-1)\delta}+ Z_t,
\end{equation}
where $\phi=e^{-\alpha\delta}$ and $Z_t$ is white noise with mean 0 and variance $\sigma^2_z=\frac{\sigma^2(1-e^{-2\alpha\delta})}{2\alpha}$.
When $\alpha$ or $\delta$ is small, the sampled discrete time first order autoregressive process approaches to near unit root solution. In such a case, large estimation bias for the parameter $\phi$ is well known. For fixed $\alpha$ and small $\delta$, \citet{TC} showed that maximum likelihood estimate of the drift parameter $\alpha$ is  $O(\frac1{n\delta})$ where $n\delta$ is referred as total span of observation which increases indefinitely as sample size increases. The bias corrected maximum likelihood estimator for the uniformly sampled data has been studied in details (see \citet{Yu}).

In some situations, the continuous time AR(1) process is also used to model irregularly observed continuous time phenomenon. Further in some application where one has controlled over the sampling mechanism, stochastic point process sampling has been used to observe continuous time process and subsequently OU model is fitted to the data \citep{SM,TAA}.
In such a case, when one can design the sampling time points to observe the underlying continuous time process, a natural constraint of the minimum separation  between successive samples often occurs due to technological or economic consideration \citep{HMR}. \citet{SS} studied the effect of this constraint on the spectrum estimation of continuous time stationary stochastic processes. Here, we study the effect of such a constraint on the estimation of the drift parameter of the continuous time AR(1) process.

\section{Estimation of parameters}\label{s2}

Let $X_t$ be the stationary solution of the diffusion equation (\ref{AReqn}) and $\tau=\{t_k; k\in \mathbb{Z}\}$ be the sequence of sampling time points. Let $\Delta_k=t_k-t_{k-1}$. When the derivative of continuous time noise is considered as Wiener process then using the normality of the noise, the conditional distribution of $X_{t_{i+1}}$ given $\{X_{t_i},t_{i+1},t_i\}$, for $i\ge 1$, is given by
\begin{eqnarray*}
X_{t_{i+1}}\bigg|\{X_{t_i},t_{i+1},t_i\}\thicksim N\left(X_{t_i}e^{-\alpha(t_{i+1}-t_i)},\frac{\sigma^2}{2\alpha}\left(1-e^{-2\alpha(t_{i+1}-t_i)}\right)\right).
\end{eqnarray*}
The log likelihood of the data is
\begin{eqnarray*}
&&\hskip-30pt L(\alpha,\sigma^2|\tau)\\&=&-\frac{n}{2}\log\frac{\sigma^2}{2\alpha} -\frac{1}{2}\sum_{k=2}^{n}\log(1-e^{-2\alpha\Delta_k})
-\frac{1}{2}\sum_{k=2}^{n}\frac{(X_{t_k}-e^{-\alpha\Delta_k}X_{t_{k-1}})^2}{\frac{\sigma^2}{2\alpha}\left(1-e^{-2\alpha\Delta_k}\right)}+c,
\end{eqnarray*}
where $c$ is a constant.

When the uniform point sampling is used to observe the process with the inter sample spacing $\Delta_k=\delta$, the maximum likelihood estimator (mle) for the parameters is  given by
\begin{eqnarray*}
\hat\alpha_{u}&=&-\frac{\log(\hat\phi)}{\delta}\\
\hat\sigma^2_{u}&=&\frac{2\hat\alpha_{u}\frac{1}{n}\sum_{i=1}^{n-1}(X_{t_{i+1}}-\hat\phi X_{t_i})^2}{1-\hat\phi^2},
\end{eqnarray*}
where $\hat\phi=\frac{\sum_{i=1}^{n-1}X_{t_{i+1}}X_{t_i}}{\sum X_{t_i}^2}$.

We now turn to the stochastic point sampling whose inter-sample spacing $\Delta$ is independent from the process $X$ and constitutes a sequence of identically, independently distributed random variable with probability density function $f(\cdot)$ (say). When such a sampling scheme is used to observe the process, the maximum likelihood estimators of the parameter can not be expressed explicitly but evaluated by numerical methods. Further, to the best of our knowledge, an optimal choice of the density as well as sampling rate is not known in general.

The maximum likelihood estimates are derived using the distributional properties of the noise. We propose a moment estimator for the drift parameter and analyze its property. The consistency of proposed moment estimator is derived by using the property of stationarity. Thus, this is a distribution free estimator for the drift of the continuous time AR(1) process.

By using the properties of independent and identically distributed sequence of inter sample spacing, a moment estimator of the drift and innovation variance of the OU process is derived as follows.
\begin{eqnarray}
E [X_{t_{i+1}}X_{t_i}]&=&E_t\left[E\{X_{t_{i+1}}X_{t_i}|t_{i+1},t_{i}\}\right]
=E_t\left[\frac{\sigma^2}{2\alpha}e^{-\alpha (t_{i+1}-t_i)}\right]
=\frac{\sigma^2}{2\alpha}\int_{0}^{\infty}e^{-\alpha t}f(t)dt\notag\\
&=&\frac{\sigma^2}{2\alpha}g(\alpha)\label{mo_d},
\end{eqnarray}
where $g(\alpha)=\int_{0}^{\infty}e^{(-\alpha t)}f(t)dt$ is the Laplace transform of the inter sample spacing variable~$\Delta$.
Note that $g(\alpha)$ is monotonic decreasing as
\begin{eqnarray*}
\frac{d}{d\alpha}g(\alpha)&=&\int_{0}^{\infty}\frac{de^{-\alpha t}}{d\alpha} f(t) dt
=-\int_{0}^{\infty}e^{-\alpha t} t f(t) dt<0.
\end{eqnarray*}
Here, the change of differential and integral operation is possible by virtue of the Dominated Convergence Theorem.
Similarly, by using
\begin{eqnarray}\label{mo_v}
E[X_{t_i}^2]&=&\frac{\sigma^2}{2\alpha},
\end{eqnarray}
we consider the sample moments
\begin{eqnarray*}
T_n&=&\frac{1}{n}\sum_{i=1}^{n-1}X_{t_{i+1}}X_{t_i};\\
V_n&=&\frac{1}{n}\sum_{i=1}^{n}X_{t_i}^2.
\end{eqnarray*}
The moment estimator of the drift and the variance parameter is given as
\begin{eqnarray}
\hat\alpha_n&=&g^{-1}\left(\frac{T_n}{V_n}\right)\\
\hat\sigma^2_n&=&2\hat\alpha_n V_n.
\end{eqnarray}

\section{Consistency of the moment estimator}
In this section, we establish the consistency of the moment estimators proposed in Section~\ref{s2}. Let $\eta=\frac{\sigma^2}{2\alpha}$.

\bigskip\noindent
{\bf Proposition 1.} {\it If the mean of inter sample spacing is finite, we have
\begin{eqnarray}
E[ T_n]&=&\left(1-\frac1n\right)\eta g(\alpha),\\
E[V_n]&=&\eta,\\
\lim_{n\rightarrow\infty}nVar[ T_n]&=& \eta^2\left[1+g(2\alpha)+4g^2(\alpha)\left(1+\int_{0}^{\infty}e^{-2\alpha v}H(v) dv\right)\right]+o\left(1 \right),\qquad\mbox{}\\
\lim_{n\rightarrow\infty}nVar[V_n]&=& \eta^2 \left[2+4\int_0^{\infty}e^{-2\alpha v}H(v)dv\right]+o\left(1 \right),\\
\lim_{n\rightarrow\infty}nCov(T_n, V_n)&=&\eta^2 \left[4g(\alpha)+4g(\alpha)\int_0^{\infty}e^{-2\alpha v}H(v)dv\right] +o\left(1 \right),
\end{eqnarray}
where $H(v)=\sum_{i=1}^{\infty}f^{(i)}(v)$ where $f^{(k)}(\cdot)$ is the the k-fold convolution of the density $f(\cdot)$.}

\bigskip\noindent
{\bf Proposition 2.} {\it If the mean of inter sample spacing is finite, we have
\begin{eqnarray*}
\lim_{n\rightarrow\infty}n\left\{E[\widehat{g}(\alpha)]-g(\alpha)\right\}&=&-3g(\alpha)+o\left(1\right)\\
\lim_{n\rightarrow\infty}n~Var[\widehat{g}(\alpha)]&=&1+g(2\alpha)-2g^2(\alpha)+o\left(1\right).
\end{eqnarray*}}

\noindent
We use Proposition~1 and~2 to obtain the rate of convergence for the moment estimator of the drift and innovation variance parameter.
 We have the following result.

\bigskip\noindent
{\bf Theorem 1} {\em Under the conditions of Proposition 1 and 2,
\begin{enumerate}
\item The bias of the estimator is given by \begin{eqnarray*}
\lim_{n\rightarrow\infty}n\left\{E[\hat\alpha]-\alpha\right\}&\!\!=\!\!&-\frac{3g(\alpha)}{g'(\alpha)}
-\frac{g''(\alpha)[1+g(2\alpha)-2g^2(\alpha)]}{2(g'(\alpha))^3}+o\left(1\right),\\
\lim_{n\rightarrow\infty}n\left\{E[\hat\sigma^2]-\sigma^2\right\}&\!\!=\!\!&2\eta\left\{-\frac{g(\alpha)}{g'(\alpha)}-\frac{g''(\alpha)[1+g(2\alpha)-2g^2(\alpha)]}{2(g'(\alpha))^3}\right\}
+o\left(1\right).
\end{eqnarray*}
\item The variance of the estimator is given by \begin{eqnarray*}
&&\hskip-30pt\lim_{n\rightarrow\infty}nVar[\hat\alpha]\\&\!\!=\!\!&\frac{1+g(2\alpha)-2g^2(\alpha)}{(g'(\alpha))^2}+o\left(1\right),\\
&&\hskip-30pt\lim_{n\rightarrow\infty}nVar[\hat\sigma^2]\\&\!\!=\!\!&4\eta^2\left\{ \!2\alpha^2\!+\!4\alpha^2\int_{0}^{\infty}e^{-2\alpha v}H(v)dv\!+\!\frac{4\alpha g(\alpha)}{g'(\alpha)}\!+\!\frac{1\!+\!g(2\alpha)\!-\!2g^2(\alpha)}{(g'(\alpha))^2}\!\right\}\!+\!o\left(1\right).
\end{eqnarray*}
\end{enumerate}}

Theorem~1 establishes the consistency of the moment estimator.

\section{Optimal sampling rate}\label{s4}
Here, we consider two special case of inter-sample spacing distribution and look for the optimal sampling rate which minimizes the bias and variance of the estimator. First, we consider inter-sample spacing distribution to be exponential. Second, under the constraint on the minimum separation between successive samples (say, $\delta$), we consider the inter-sample spacing distribution to be truncated exponential.

When the inter sample spacing is exponential distribution with rate $\beta$, we have
$$g(\alpha)=\frac{\beta}{\beta+\alpha}.$$
Then by using Theorem~1, we have
\begin{eqnarray*}
\lim_{n\rightarrow\infty}n Bias(\hat \alpha)&=&\frac{(\beta+\alpha)(2\alpha^3+3\beta^3+8\alpha\beta^2+6\alpha^2\beta)}{\beta^2(\beta+2\alpha)}\\
\lim_{n\rightarrow\infty}n Var(\hat\alpha)&=&\frac{2\alpha(\beta+\alpha)^2(\beta^2+\alpha^2+3\alpha\beta)}{\beta^2(\beta+2\alpha)}.
\end{eqnarray*}
It is evident form the expression of the asymptotic bias and variance of the estimator that the different average sampling rate $\beta$ minimizes the respective expression. However, the optimal rate corresponding to the asymptotic bias and variance appears not to vary much. Figure~1 shows the plot of optimal sampling rate corresponding to the asymptotic bias and variance. The problem of large bias in the drift estimation, when drift is small and process is observed at uniformly spaced time point, is well known. If it is known a priory that the drift is contained in a interval where the upper end of the interval is small, we can choose the optimal sampling rate which minimizes the maximum relative bias over the interest of the drift interval.
\begin{figure}
\begin{center}
 \includegraphics[width=4in]{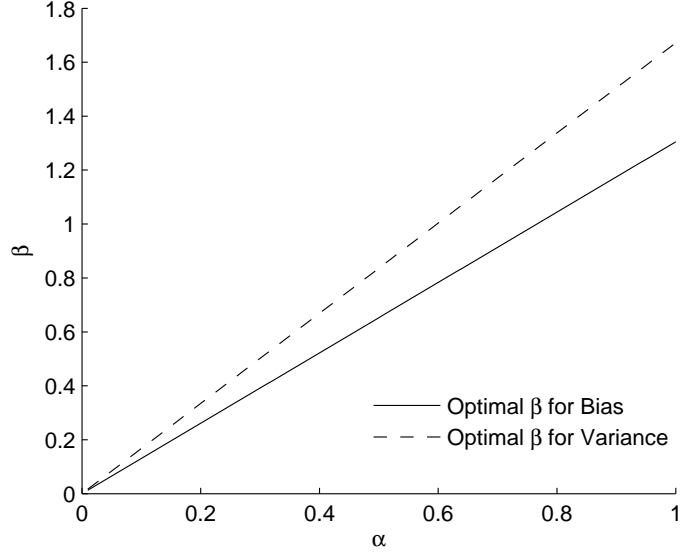}
 \caption{\small Optimal average sampling rate corresponding to bias and variance for Poisson}
\end{center}
\end{figure}

We now turn to the effect of constraint on the minimum separation between successive samples. Let the inter-sample spacing be $\delta+ \Delta$ where $\delta>0$ is the minimum threshold between successive samples and $\Delta$ is exponential with rate $\beta$. Then, we have
$$g(\alpha)=\frac{\beta e^{-\alpha\delta}}{\beta+\alpha}.$$
The expression for the asymptotic bias and variance can be obtained from Theorem~1. However, these expression turns out to be complicated. Figure~2 shows the plot of optimal sampling rate corresponding to the asymptotic bias and variance when the inter sample spacing is exponential truncated from left at $\delta=0.1$ and $0.5$.
\begin{figure}
\begin{center}
 \includegraphics[width=4in]{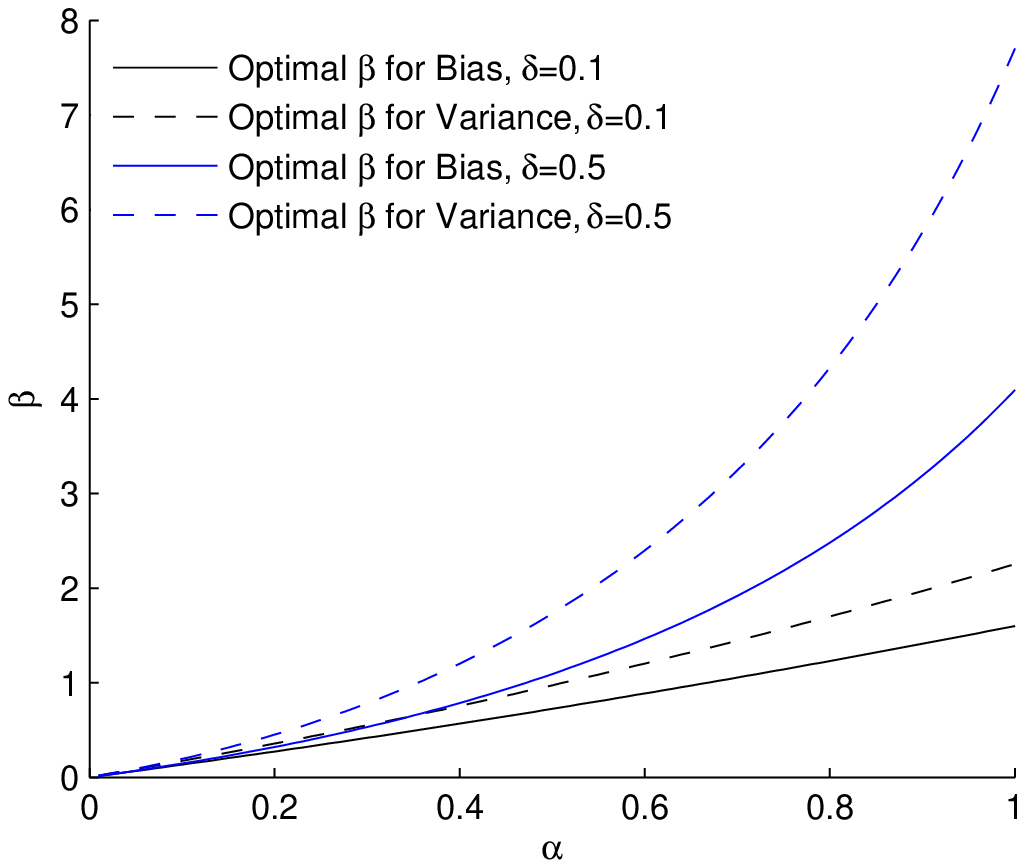}
 \caption{\small Optimal average sampling rate corresponding to bias and variance for truncated Poison}
\end{center}
\end{figure}

Figure~2 indicates that the difference between the optimal $\beta$ corresponding to asymptotic bias and variance appears to be large as the minimum separation between successive samples $\delta$ increases. This also indicates that the relative bias of the drift is increasing function of the minimum separation $\delta$.

\section{Conclusion}
The article proposes a distribution free moment estimator for the drift and the innovation variance of the continuous time AR(1) process. The expression for the asymptotic bias and variance of the proposed moment estimators of the parameters are derived. The problem of large bias in the estimation of the drift parameter is well known when drift is small. When the interest of the drift interval is known a priory, we can use these expression to choose the optimal average sampling rate while designing the stochastic sampling time points as described in the Section~\ref{s4}. The constraint on the minimum separation between the successive samples influences the choice of optimal average sampling rate. It appears when the minimum separation between successive samples is large the relative asymptotic bias is large in particular when drift is small.

\section{Appendix}
{\bf Proof of Proposition 1}\\
{\bf \it Part (i)}. From (\ref{mo_d}), we have
$$E(T_n)=\frac{1}{n}\sum_{i=1}^{n-1}E\left[X_{t_{i+1}}X_{t_i}\right]=\left(1-\frac1n\right)\eta g(\alpha).$$
\hfill\(\Box\)

\noindent
{\bf \it Part (ii)}. From (\ref{mo_v}), we have
$$E(V_n)=\frac{1}{n}\sum_{i=1}^{n}E\left[X_{t_i}^2\right]=\eta.$$
\hfill\(\Box\)

\noindent
{\bf \it Part (iii)}.
By using the Normality of the process $X$, we have
\begin{eqnarray}
&&\hskip-20pt Var(T_n)\notag\\&\!\!=\!\!&E_t\left[\frac{1}{n^2}\sum_{i=1}^{n-1}\sum_{j=1}^{n-1}Cov(X_{t_{i+1}},X_{t_j} )Cov(X_{t_{j+1}},X_{t_i} )+Cov(X_{t_{i+1}},X_{t_{j+1}} )Cov(X_{t_i},X_{t_j} )\right]\notag\\
 &\!\!=\!\!&\frac{1}{n^2}\sum_{i=1}^{n-1}\sum_{j=1}^{n-1} \left(\frac{\sigma^2}{2\alpha}\right)^2 \left\{E_t\left[e^{-\alpha|t_{i+1}-t_j|}e^{-\alpha|t_{j+1}-t_i|}\right]+E_t\left[e^{-\alpha|t_{i+1}-t_{j+1}|}e^{-\alpha|t_{j}-t_i|}\right]\right\}\notag\\
 &\!\!=\!\!&I_1+I_2, (\mbox{\it say}).
\end{eqnarray}
We first consider the term $I_1$. Note that
\begin{eqnarray*}
 I_1&=&\frac{\eta^2}{n^2}\sum_{i=2}^{n-1}\sum_{j=1}^{i-1}E_t\left[e^{-\alpha(t_{i+1}-t_j)}e^{-\alpha(t_i-t_{j+1})}\right]
+\frac{\eta^2}{n^2}\sum_{i=1}^{n-1}E_t\left[e^{-2\alpha(t_{i+1}-t_i)}\right]\\
&&+\frac{\eta^2}{n^2}\sum_{i=1}^{n-1}\sum_{j= i+1}^{n-2}E_t\left[e^{-\alpha(t_j-t_{i+1})}e^{-\alpha(t_{j+1}-t_i)}\right]\\
&=&I_{11}+I_{12}+I_{13}, (\mbox{say}).
\end{eqnarray*}
Note that
\begin{eqnarray*}
I_{11}&\!\!=\!\!&\frac{\eta^2}{n^2}\left\{\sum_{i=3}^{n-1}\sum_{j=1}^{i-2}E_t\left[e^{-\alpha(t_{i+1}-t_j)}e^{-\alpha(t_i-t_{j+1})}\right]
+\sum_{i=2}^{n-1}E_t\left[e^{-\alpha(t_{i+1}-t_{i-1})}\right]\right\}\\
&\!\!=\!\!&\frac{\eta^2}{n^2}\left\{\sum_{i=3}^{n-1}\sum_{j=1}^{i-2}E_t\left[e^{-\alpha(t_{i+1}-t_j)}e^{-\alpha(t_i-t_{j+1})}\right]
+\sum_{i=2}^{n-1}E_t\left[e^{-\alpha(t_{i+1}-t_i+t_i-t_{i-1})}\right]\right\}.
\end{eqnarray*}
Let  $\vartheta_1=t_{i+1}-t_i$, $\vartheta_2=t_i-t_{j+1}$ and $\vartheta_3=t_{j+1}-t_j$. By using the independence of the
inter sample spacing, the probability density function of $\vartheta_i$ are $f(\cdot)$, $f^{(i-j-1)}(\cdot)$ and $f(\cdot)$ respectively, where $f^{(l)}(\cdot)$
is the $l$ fold convolution of the density $f(\cdot)$.
Thus, we have
\begin{eqnarray*}
 I_{11}&\!\!\!\!=\!\!\!\!&\frac{\eta^2}{n^2}\sum_{i=3}^{n-1}\sum_{j=1}^{i-2}
\int_{0}^{\infty}\int_{0}^{\infty}\int_{0}^{\infty}e^{-\alpha(v_1+v_2+v_3)}e^{-\alpha(v_2)} f(v_1)f^{(i-j-1)}(v_2)f(v_3)dv_1dv_2dv_3\\
&&+\frac{\eta^2}{n^2}\sum_{i=2}^{n-1}\int_{0}^{\infty}\int_{0}^{\infty}e^{-\alpha(v_1+v_3)}f(v_1)f(v_3)dv_1dv_3\\
&\!\!\!\!=\!\!\!\!&\frac{\eta^2}{n}\!\int_{0}^{\infty}\!\int_{0}^{\infty}\!\int_{0}^{\infty}\!e^{-\alpha(v_1\!+\!v_2\!+\!v_3)}e^{-\alpha(v_2)} f(v_1)\!\left(\!\frac{1}{n}\sum_{i=3}^{n-1}\sum_{j= 1}^{i-2}f^{(i\!-\!j\!-\!1)}(v_2)\!\right)\!f(v_3)dv_1dv_2dv_3\\
&&+\frac{\eta^2}{n}\left(1-\frac2n\right)g^2(\alpha).
\end{eqnarray*}
Note that
\begin{eqnarray*}
\frac{1}{n}\sum_{i=3}^{n-1}\sum_{j=1}^{i-2}f^{(i-j-1)}(v_2)&=&\frac{1}{n}\sum_{j=1}^{n-3}\sum_{l=1}^{n-j-2}f^{(l)}(v_2)=\frac{1}{n}\sum_{l=1}^{n-3}\sum_{j=1}^{n-l-2}f^{(l)}(v_2)\\
&=& \sum_{l=1}^{n-3} \left(1-\frac{l+2}{n}\right)f^{(l)}(v_2).
\end{eqnarray*}
Thus we have
\begin{eqnarray*}
nI_{11}&\!\!\!=\!\!\!&\eta^2 g^2(\alpha) \!\left\{1\!+\!\int_{0}^{\infty}e^{-\alpha v_2}\sum_{l=1}^{n-3} f^{(l)}(v_2)dv_2\!-\!\frac1n \int_{0}^{\infty}e^{-\alpha v_2}\sum_{l=1}^{n-3} (l\!+\!2)f^{(l)}(v_2)dv_2\!-\!\frac2n\right\}.
\end{eqnarray*}
Let $H(v)=\sum_{l=1}^{\infty}f^{(l)}(v)$. Note that $H(\cdot)$ is bounded as the mean of the inter sample spacing is finite. Further,
$$\sum_{l=1}^{\infty}f^{(l)}(v)=\frac12 +\int_{0}^{v} H(t)dt-\frac12 H(v).$$
By using Dominated Convergence Theorem (DCT), we have
\begin{eqnarray*}
\lim_{n\rightarrow\infty}n I_{11}&=&\eta^2g^2(\alpha)\left(1+\int_{0}^{\infty}e^{-2\alpha v}H(v) dv\right)+o\left(\frac1n\right).
\end{eqnarray*}
Now we turn to $I_{12}$. Let $t_{i+1}-t_i=\vartheta_1$, then
\begin{eqnarray*}
I_{12}&=&\eta^2\frac{1}{n}\left(1-\frac1n\right)\int_{0}^{\infty}e^{-2\alpha}f(v)dv=\frac{1}{n}\left(1-\frac1n\right)\eta^2g(2\alpha).
\end{eqnarray*}
Thus, we have
$$\lim_{n\rightarrow\infty}n I_{12}=\eta^2g(2\alpha)+o\left(\frac1n\right). $$
Note that the term $I_{13}$ is a symmetric to $I_{11}$. By using this symmetry, we have
\begin{eqnarray*}
\lim_{n\rightarrow\infty}n I_{13}&=&\eta^2g^2(\alpha)\left(1+\int_{0}^{\infty}e^{-2\alpha v}H(v) dv\right)+o\left(\frac1n\right).
\end{eqnarray*}
Thus,
$$\lim_{n\rightarrow\infty}nI_1=\eta^2\left[g(2\alpha)+2g^2(\alpha)\left(1+\int_{0}^{\infty}e^{-2\alpha v}H(v) dv\right)\right]+o\left(\frac1n\right).$$
We now consider $I_2$. Note that
\begin{eqnarray*}
I_2&\!\!\!=\!\!\!&\frac{\eta^2}{n}\left(1-\frac1n\right)+\frac{\eta^2}{n^2}
\sum_{i=2}^{n-1}\sum_{j=1}^{i-1}E_t\left[e^{-\alpha(t_{i+1}-t_{j+1})}e^{-\alpha(t_i-t_{j})}\right]\\
&&+\frac{\eta^2}{n^2}\sum_{i=1}^{n-1}\sum_{j= i+1}^{n-2}E_t\left[e^{-\alpha(t_{j+1}-t_{i+1})}e^{-\alpha(t_{j}-t_i)}\right]\\
&\!\!\!=\!\!\!&I_{21}+I_{22}+I_{23}.
\end{eqnarray*}
Consider $I_{22}$ and let $\vartheta_1=t_{i+1}-t_i$, $\vartheta_2=t_i-t_{j+1}$ and $\vartheta_3=t_{j+1}-t_j$, then
\begin{eqnarray*}
I_{22}&\!\!\!\!=\!\!\!\!&\frac{\eta^2}{n^2}\sum_{i=3}^{n-1}\sum_{j=1}^{i-2}
 \int_{0}^{\infty}\!\int_{0}^{\infty}\!\int_{0}^{\infty}\!e^{-\alpha(v_1\!+\!v_2)}e^{-\alpha(v_2\!+\!v_3)} f(v_1)f^{(i\!-\!j\!-\!1)}(v_2) f(v_3)dv_1dv_2dv_3\\
&&+\frac{\eta^2}{n^2}\sum_{i=2}^{n-1}\int_{0}^{\infty}\!\int_{0}^{\infty}\!e^{-\alpha v_1}e^{-\alpha v_3}f(v_1)f(v_3)dv_1dv_3\\
&\!\!\!\!=\!\!\!\!&\frac{\eta^2}{n^2} g^2(\alpha)
 \int_{0}^{\infty}e^{-2\alpha v_2}\sum_{l=1}^{n-3} \left(1\!-\!\frac{l\!+\!2}{n}\right)f^{(l)}(v_2)dv_2\!+\!\frac{\eta^2}n\left(1\!-\!\frac2n\right)g^2(\alpha).
\end{eqnarray*}
By using a similar argument as in case of $I_{12}$, we have
$$\lim_{n\rightarrow\infty}nI_{22}=\eta^2g^2(\alpha)\left(1+\int_{0}^{\infty}e^{-2\alpha v}H(v) dv\right)+o\left(\frac1n\right).$$
By using symmetry of $I_{22}$ and $I_{23}$, we have
$$\lim_{n\rightarrow\infty}nI_{23}=\eta^2g^2(\alpha)\left(1+\int_{0}^{\infty}e^{-2\alpha v}H(v) dv\right)+o\left(\frac1n\right).$$
By combining all the terms, we have
\begin{eqnarray*}
\lim_{n\rightarrow\infty} nVar[T_n]&=&\eta^2\left[1+g(2\alpha)+4g^2(\alpha)\left(1+\int_{0}^{\infty}e^{-2\alpha v}H(v) dv\right)\right]+o\left(\frac1n\right).
\end{eqnarray*}
\hfill\(\Box\)

\bigskip\noindent
{\bf \it Part (iv)}. Note that
\begin{eqnarray*}
Var(V_n)&=&\frac{2}{n^2}\sum_{i=1}^{n}\sum_{j=1}^{n}E_t\left[Cov(X_{t_1},X_{t_j})^2\right]=\frac{2\eta^2} {n^2}\sum_{i=1}^n\sum_{j=1}^nE_t[e^{-2\alpha|t_i-t_j|}]\\
&=&\frac{2\eta^2} {n}+\frac{2\eta^2} {n^2}\sum_{i=2}^n\sum_{j=1}^{i-1}E_t[e^{-2\alpha(t_i-t_j)}]
+\frac{2\eta^2} {n^2}\sum_{i=1}^n\sum_{j=i+1}^{n-1}E_t[e^{-2\alpha(t_j-t_i)}]\\
&=&J_1+J_2+J_3.\\
\end{eqnarray*}
Consider $J_2$ and let $\vartheta_1=t_i-t_j$, then
\begin{eqnarray*}
J_2=\frac{2\eta^2} {n^2}\sum_{i=2}^n\sum_{j=1}^{i-1}\int_0^\infty e^{-2\alpha v}f^{(i-j)}(v)dv
=\frac{2\eta^2} {n}\int_0^\infty e^{-2\alpha v}\left\{\sum_{l=1}^{n-1}\left(1-\frac ln\right)f^{(l)}(v)\right\}dv.
\end{eqnarray*}
By using a similar argument as in case of $I_{12}$, we have
$$\lim_{n\rightarrow\infty}nJ_2= 2 \eta^2\int_0^{\infty}e^{-2\alpha v}H(v)dv+o\left(\frac 1n\right).$$
By using the symmetry of $J_2$ and $J_3$, we have
$$\lim_{n\rightarrow\infty}nJ_3= 2\eta^2\int_0^{\infty}e^{-2\alpha v}H(v)dv+o\left(\frac 1n\right).$$
By combining the terms, we have
\begin{eqnarray*}
\lim_{n\rightarrow\infty} nVar[V_n]&=&\eta^2\left[2+4\int_{0}^{\infty}e^{-2\alpha v}H(v) dv\right]+o\left(\frac1n\right).
\end{eqnarray*}
\hfill\(\Box\)

\bigskip\noindent
{\bf \it Part (v)}. Note that
\begin{eqnarray*}
&&\hskip-20pt Cov(T_n,V_n)
\\&=&\frac{2\eta^2}{n^2}\sum_{i=1}^{n-1}\sum_{j=1}^{n-1}E_t\left[e^{-\alpha|t_{i+1}-t_j|}e^{-\alpha|t_{i}-t_j|}\right]
+\frac{2\eta^2}{n^2}\sum_{i=1}^{n-1}E_t\left[e^{-\alpha(t_n-t_{i+1})}e^{-\alpha(t_{n}-t_i)}\right]\\
&=&\frac{2\eta^2}{n^2}\sum_{i=1}^{n-1}E_t\left[e^{-\alpha (t_{i+1}-t_i)}\right]
+\frac{2\eta^2}{n^2}\sum_{i=2}^{n-1}\sum_{j=1}^{i-1}E_t\left[e^{-\alpha(t_{i+1}-t_j)}e^{-\alpha(t_{i}-t_j)}\right]\\
&&+\frac{2\eta^2}{n^2}\sum_{i=1}^{n-2}\sum_{j=i+1}^{n-1}E_t\left[e^{-\alpha(|t_j-t_{i+1}|)}e^{-\alpha(t_{j}-t_i)}\right]
+\frac{2\eta^2}{n^2}\sum_{i=1}^{n-1}E_t\left[e^{-\alpha(t_n-t_{i+1})}e^{-\alpha(t_{n}-t_i)}\right]\\
&=&J_1+J_2+J_3+J_4.
\end{eqnarray*}
Consider $J_1$ and let $t_{i+1}-t_i=\vartheta_1$, then
\begin{eqnarray*}
J_1&=&\frac{2(n-1)}{n^2}\eta^2 \int_0^\infty e^{-\alpha v}f(v)dv
=\frac{2}{n}\left(1-\frac1n\right)\eta^2g(\alpha).
\end{eqnarray*}
Now consider $J_2$ and let $t_{i+1}-t_i=\vartheta_1$ and $t_{i}-t_j=\vartheta_2$, then
\begin{eqnarray*}
J_2&=&\frac{2\eta^2}{n^2}\sum_{i=2}^{n-1}\sum_{j=1}^{i-1}\int_{0}^{\infty}\int_{0}^\infty e^{-\alpha(v_1+v_2)}e^{-\alpha v_2} f(v_1)f^{(i-j)}(v_2) dv_1dv_2\\
&=&\frac{2\eta^2}{n}g(\alpha)\int_{0}^\infty e^{-2\alpha v_2} \sum_{l=1}^{n-2}\left(1-\frac{2+l}{n}\right)f^{(l)}(v_2) dv_2.
\end{eqnarray*}
By using a similar arguments as in case of $I_{12}$, we have
\begin{eqnarray*}
\lim_{n\rightarrow\infty}nJ_2=2\eta^2 g(\alpha)\int_0^{\infty}e^{-2\alpha v}H(v)dv+o\left(\frac1n\right).
\end{eqnarray*}
Consider $J_3$ and let $t_j-t_{i+1}=\vartheta_1$ and $t_{i+1}-t_i=\vartheta_2$, then
\begin{eqnarray*}
J_3&=&\frac{2\eta^2}{n^2}\sum_{i=1}^{n-3}\sum_{j=i+2}^{n-1}E_t\left[e^{-\alpha(t_j-t_{i+1})}e^{-\alpha(t_{j}-t_i)}\right]+\frac{2\eta^2}{n^2}
\sum_{i=1}^{n-2}E_t\left[e^{-\alpha(t_{i+1}-t_i)}\right]\\
&=&\frac{2\eta^2}{n^2}\sum_{i=1}^{n-3}\sum_{j=i+2}^{n-1}\int_0^\infty\int_0^\infty e^{-\alpha v_1}e^{-\alpha(v_1+v_2)} f^{(j-i-1)}(v_1)f(v_2)dv_1dv_2\\
&&+\frac{2(n-2)}{n^2}\eta^2\int_0^\infty e^{-\alpha v_2}f(v_2) dv_2\\
&=&\frac{2\eta^2}{n}g(\alpha)\int_0^\infty e^{-2\alpha v_1}\sum_{l=1}^{n-3}\left(1-\frac{2+l}{n}\right)f^{(l)}(v_1)dv_1+\frac{2}{n}\left(1-\frac2n\right)\eta^2g(\alpha).
\end{eqnarray*}
By using a similar arguments as in case of $I_{12}$, we have
\begin{eqnarray*}
\lim_{n\rightarrow\infty}nJ_3&=&\eta^2 \left[2g(\alpha)\int_0^{\infty}e^{-2\alpha v}H(v)dv+2g(\alpha)\right]+o\left(\frac1n\right).
\end{eqnarray*}
Now consider $J_4$ and let $t_{i+1}-t_i=\vartheta_1$ and $t_n-t_{i+1}=\vartheta_2$, then
\begin{eqnarray*}
J_4&\!\!\!=\!\!\!&\frac{2\eta^2}{n^2}\sum_{i=1}^{n-2}\int_{0}^{\infty}\int_{0}^\infty e^{-\alpha v_2}e^{-\alpha (v_2+v_1)}f^{(n-i-1)}(v_2)f(v1)dv_1dv_2
+\frac{2\eta^2}{n^2}\int_{0}^{\infty}e^{-\alpha v_2}f(v_2)dv_2\\
&\!\!\!=\!\!\!&\frac{2\eta^2}{n^2} g(\alpha)\left[\int_{0}^\infty e^{-2\alpha v_2}
\sum_{l=1}^{n-2}f^{(l)}(v_2)dv_2+ 1\right]=O\left(\frac1{n^2}\right).
\end{eqnarray*}
By combining these terms, we have
\begin{eqnarray*}
\lim_{n\rightarrow\infty}nCov(T_n,V_n)&=&\eta^2 \left[4g(\alpha)+4g(\alpha)\int_0^{\infty}e^{-2\alpha v}H(v)dv\right]+o\left(\frac1n\right).
\end{eqnarray*}
This completes the proof.
\hfill\(\Box\)

\bigskip\noindent
{\bf Proof of Proposition 2}.
Note that by using the first order approximation, $\hat g(\alpha)=\frac{T_n}{V_n}$ can be expressed as
\begin{eqnarray}\label{g_hat_app}
\widehat{g}(\alpha)&\!\!=\!\!&g(\alpha)\!+\!\frac{1}{\eta}(T_n\!-\!\eta g(\alpha))\!-\!\frac{ g(\alpha)}{\eta}(V_n\!-\!\eta)+\frac{ g(\alpha)}{\eta^2}(V_n\!-\!\eta)^2\!-\!\frac{1}{\eta^2}(T_n\!-\!\eta g(\alpha))(V_n\!-\!\eta).\notag\\
\end{eqnarray}
From Proposition~1 and (\ref{g_hat_app}), we have
\begin{eqnarray*}
E\left[\widehat{g}(\alpha)\right]-g(\alpha)&=&\frac{1}{\eta}E(T_n-\eta g(\alpha))+\frac{ g(\alpha)}{\eta^2}Var(V_n)-\frac{1}{\eta^2}Cov(T_n,V_n)\\
&=&-\frac{3g(\alpha)}{n}+o\left(\frac1n\right).
\end{eqnarray*}
Similarly by using Proposition~1 and (\ref{g_hat_app}), we have
\begin{eqnarray*}
Var\left[\widehat{g}(\alpha)\right]&=&\frac{1}{\eta^2}Var(T_n)+\left(\frac{ g(\alpha)}{\eta}\right)^2Var(V_n)-2\frac{g(\alpha)}{\eta^2}Cov(T_n,V_n)\\
&=&\frac{1+g(2\alpha)-2g^2(\alpha)}n+o\left(\frac1n\right).
\end{eqnarray*}
\hfill\(\Box\)

\bigskip\noindent
{\bf Proof of Theorem 1.}
Note that  $\hat\alpha=g^{-1}\left(\frac{T_n}{V_n}\right)=g^{-1}\left(\widehat g(\alpha)\right)$.
Further, we have
\begin{eqnarray*}
\frac{dg^{-1}(x)}{dx}=\frac{1}{g'\left(g^{-1}(x)\right)}; && \frac{d^2g^{-1}(x)}{dx^2}=-\frac{g''\left(g^{-1}(x)\right)}{\left[g'\left(g^{-1}(x)\right)\right]^3},
\end{eqnarray*}
where $g'(\cdot)$ and $g''(\cdot)$ are the first and second order derivative of $g$. By using Taylor Series approximation, we have
\begin{eqnarray}
\hat\alpha&=&\alpha+\frac{1}{g'(\alpha)}(\widehat{g}(\alpha)-g(\alpha))-\frac{g''(\alpha)}{2(g'(\alpha))^3}(\widehat{g}(\alpha)-g(\alpha))^2\label{alpha_hat_app}\\
\hat\sigma^2&=&2\hat\alpha V_n=
\sigma^2+2\alpha(V_n-\eta)+2\eta(\hat\alpha-\alpha)+2(V_n-\eta)(\hat\alpha-\alpha)\label{sigma_hat_app}.
\end{eqnarray}
{\bf \it Part (i)}. By using Proposition 2 and (\ref{alpha_hat_app}), we have
\begin{eqnarray*}
E[\hat\alpha]-\alpha&=&\frac{1}{g'(\alpha)}E(\widehat{g}(\alpha)-g(\alpha))-\frac{g''(\alpha)}{2(g'(\alpha))^3}Var(\widehat{g}(\alpha))\\
&=&\frac{1}{n}\left\{-\frac{3g(\alpha)}{g'(\alpha)}
-\frac{g''(\alpha)[1+g(2\alpha)-2g^2(\alpha)]}{2(g'(\alpha))^3}\right\}+o\left(\frac1n\right).
\end{eqnarray*}
Note that by using (\ref{sigma_hat_app}), we have
\begin{eqnarray}\label{bias_sigma_hat}
E[\hat\sigma^2]-\sigma^2&=&2\eta E(\hat\alpha-\alpha)+2Cov(V_n,\hat\alpha).
\end{eqnarray}
Further, by using the following first order approximation, we have
\begin{eqnarray*}
\hat\alpha-\alpha&=&\frac{1}{g'(\alpha)}(\widehat{g}(\alpha)-g(\alpha))=\frac{1}{g'(\alpha)\eta}(T_n-\eta g(\alpha))-\frac{g(\alpha)}{g'(\alpha)\eta}(V_n-\eta).
\end{eqnarray*}
By using Proposition~1, we have
\begin{eqnarray}\label{cov_a_V}
Cov(\hat\alpha,V_n)&=&\frac{1}{g'(\alpha)\eta}Cov(T_n,V_n)-\frac{g(\alpha)}{g'(\alpha)\eta}Var(V_n)=\frac{1}{n}\frac{2\eta g(\alpha)}{g'(\alpha)}+o\left(\frac{1}{n}\right).
\end{eqnarray}
Thus, by using Propositon~2 and (\ref{bias_sigma_hat}) and (\ref{cov_a_V}), we have
\begin{eqnarray*}
E[\hat\sigma^2]-\sigma^2&=&\frac{2\eta}{n}\left\{-\frac{g(\alpha)}{g'(\alpha)}-\frac{g''(\alpha)[1+g(2\alpha)-2g^2(\alpha)]}{2(g'(\alpha))^3}\right\}
+o\left(\frac1n\right).
\end{eqnarray*}
\hfill\(\Box\)

\bigskip\noindent
{\bf \it Part (ii)}. The variance expression can be easily obtained from the Proposition~1 and~2 and first order approximation made in proof of Part (i).
\hfill\(\Box\)

\bibliographystyle{apalike}
\bibliography{refrence}

\end{document}